\documentclass{amsart}
\usepackage{graphicx}

\usepackage[english]{babel}
\usepackage[cp1251]{inputenc}
\usepackage{amssymb}
\usepackage{lscape}
\usepackage{amsthm}
\usepackage{inputenc}

\newtheorem{defn}{Definition}[section]
\newtheorem{thm}{Theorem}[section]
\newtheorem{pr}{Proposition}[section]
\newtheorem{exam}{Example}[section]
\newtheorem{cor}{Corollary}[section]

\newenvironment{dem}{\rm \trivlist \item[\hskip \labelsep{\it
      Proof}:]}{\par\nopagebreak \hfill $\Box$ \endtrivlist}

\def\ll{\mathcal{L}}

\def\gg{{\mathfrak{g}}}

\begin{document}
\date{}
\author{L.M. Camacho, E.M. Ca\~{n}ete, J.R. G\'{o}mez, B.A. Omirov}

\title{\bf Quasi-filiform Leibniz algebras of maximum length.}

\maketitle
\begin{abstract}
The $n$-dimensional $p$-filiform Leibniz algebras of maximum
length have already been studied with $0\leq p\leq 2$. For Lie algebras whose nilindex is equal to $n-2$
there is only one characteristic sequence, $(n-2,1,1)$, while in Leibniz theory we obtain two possibilities: $(n-2,1,1)$ and $(n-2,2)$.
The first case (the 2-filiform case) is already known. The present paper deals with the second case, i.e., quasi-filiform non Lie Leibniz
algebras of maximum length. Therefore this work completes the study of maximum length of  Leibniz algebras with nilindex $n-p$ with $0 \leq p \leq 2$.
\end{abstract}



\section{Introduction}
The notion of length of a Lie algebra was introduced by G\'{o}mez, Jim\'{e}nez-Merch\'{a}n and Reyes in \cite{Reyes1}.
In this work they distinguished a very interesting family: algebras admitting a graduation with the greatest possible number of non-zero subspaces,
the so-called algebras of maximum length.

 Leibniz algebras appear as a natural generalization of Lie algebras (Loday, \cite{loday1} and \cite{loday}) and the concept of length can be
 defined in a similar way in this setting. It is  therefore expected that Leibniz algebras of maximum length will play a similar role to the Lie case. The cohomological properties of Leibniz algebras have been widely studied
 (see for example \cite{Dz1}- \cite{Cohomology} and \cite{ Ve}).
 A remarkable fact of the algebras of maximum length is the relative simplicity of the study of these properties \cite{LAA}.

  The study of the classification of non associative nilpotent Lie algebras is too complex.
  In fact, it appeared two centuries ago and it still remains unsolved. As to Leibniz algebras the problem
   is analogous thus we will restrict our attention to two important families of Leibniz algebras: $p$-filiform and quasi-filiform
   (see Definitions \ref{def:p-fili} and \ref{def:quasi}).

   The classification of filiform and 2-filiform Lie algebras of maximum length is given in \cite{Reyes1} and \cite{ Reyes2}
    (note that in this case there are no null-filiform algebras and the concepts of 2-filiform and quasi-filiform agree).
    In a Leibniz setting, the null-filiform case was studied in  \cite{Omirov1}, whereas the cases filiform and 2-filiform
    are developed in \cite{J.Lie.Theory2}.

   We will focus our attention on quasi-filiform Leibniz algebras.
   Let $\ll$ be an $n$-dimensional quasi-filiform non Lie Leibniz algebra,
   then its characteristic sequence is $(n-2,1,1)$ or $(n-2,2)$. The first case
   (the 2-filiform case) has been dealt with in \cite{J.Lie.Theory2}. In this work we consider the second case, i.e.,
    algebras with characteristic sequence $(n-2,2)$. Our main goal is to complete the study of algebras of maximum length with nilindex $n-2$.
     For the study of these algebras we extend the naturally graded quasi-filiform Leibniz algebras and we prove that if we consider the extension of naturally graded quasi-filiform Lie algebras then we obtain Lie algebras.

    Our main results (see Theorems \ref{th:tipo I} and
   \ref{th:tipo II} for the precise statement) will be proved in Section 3 and the following definitons will be used throughtout the present paper:

   \begin{defn}\label{def:Leibniz}
  An algebra $\ll$ over a field $F$ is called Leibniz algebra if it verifies the Leibniz identity: $[x,[y,z]]=[[x,y],z]-[[x,z],y]$ for any
  elements $x,y,z \in \ll$ and where $[,]$ is the multiplication in $\ll$.
  \end{defn}

  Note that if in $\ll$ the identity $[x,x]=0$ holds, then the Leibniz identity coincides with the
  Jacobi identity. Thus, Leibniz algebras are a generalization of Lie algebras.\\
  For a given Leibniz algebra $\ll$ we defined the following sequence: $\ll^1=\ll$ and $\ll^{k+1}=[\ll^k,\ll]$. From now on
   we will consider the complex number field. Let $\ll$ be a nilpotent algebra for which the index of nilpotency is equal to $k$.
   Let us define the naturally graded algebras as follows:
   \begin{defn}
   Let us take $\ll_i=\ll^i/\ll^{i+1}$, $1\leq i \leq n-k-1$ and $gr \ll=\ll_1 \oplus \ll_2 \oplus \dots \oplus \ll_k$. Then $[\ll_i,\ll_j]\subseteq \ll_{i+j}$
    and we obtain the graded algebra $gr \ll$. If $gr \ll$ and $\ll$ are isomorphic, in notation $gr\ll \cong \ll$, we say that $\ll$
    is a naturally graded algebra.
   \end{defn}

    The above constructed graduation is called \emph{natural graduation}.

     The set $R(\ll)=\{x \in \ll: [y,x]=0, \ \forall y \in \ll\}$ is called \emph{the right annihilator of $\ll$}. Note that for any $x, y \in \ll$ the elements $[x,x]$ and $[x,y]+ [y,x]$ are in $R(\ll)$.

  The set $Cent(\ll)=\{z \in \ll: [x,z]=[z,x]=0, \ \forall x \in \ll\}$ is called \emph{the center of $\ll$}. Note that $R(\ll)$ is an ideal of $\ll$.
  
  We define the set $\mathcal{I(\ll)}=<[x,x]: \ \forall x\in \ll>$. Note that $\mathcal{I(\ll)}$ is an ideal of $\ll$.

  Let $x$ be a nilpotent element of the set $\ll \setminus \ll^2$. For the nilpotent operator of rigth multiplication $R_x$ we define a decreasing sequence $C(x)=(n_1,n_2, \dots, n_k)$, which consists of the dimensions of Jordan
blocks of the operator $R_x$. In the set of such sequences we consider the lexicographic order, that is,
  $C(x)=(n_1,n_2, \dots, n_k)\leq C(y)=(m_1, m_2, \dots, m_s)\Longleftrightarrow$ there exists $i \in \mathbb{N}$ such that $n_j=m_j$ for any $j<i$ and $n_i<m_i$.

  \begin{defn}\label{def:char.seq}
  The sequence $C(\ll)=max C(x)_{x \in \ll \setminus \ll^2}$ is called characteristic sequence of the algebra $\ll$.
  \end{defn}

  \begin{exam}
  If $C(\ll)=(1,1, \dots,1)$, then evidently, the algebra $\ll$ is abelian.
  \end{exam}

  Let $\ll$ be an $n$-dimensional nilpotent Leibniz algebra and $p$ a non negative integer ($p<n$).
  \begin{defn}\label{def:p-fili}
  The Leibniz algebra $\ll$ is called $p$-filiform if $C(\ll)=(n-p,\underbrace{1,\dots,1}_{p})$. If $p=1$, $\ll$ is called  filiform algebra
  and if $p=0$ null-filiform algebra.
  \end{defn}

  \begin{defn}\label{def:quasi}
  A Leibniz algebra $\ll$ is called quasi-filiform if its nilindex is equal to $n-2$, namely, $\ll^{n-2} \neq \{ 0\}$ and $\ll^{n-1}=\{0\}$,
  where $n=dim(\ll)$.
  \end{defn}

    A Leibniz algebra $\ll$ is $\mathbb{Z}$-graded if $\ll=\oplus_{i \in \mathbb{Z}}V_i$,
   where $[V_i,V_j]\subseteq V_{i+j}$ for any $i,j \in \mathbb{Z}$ with a finite number of non null spaces $V_i$.

  We will say that a nilpotent Leibniz algebra $\ll $ admits the \emph{connected graduation} $\ll=V_{k_1}\oplus \dots \oplus V_{k_t}$ if $V_{k_i}\neq 0$ for any $i$ $(1 \leq i \leq t)$.
  \begin{defn}\label{def:length}
  The number $l( \oplus \ll)=l(V_{k_1}\oplus \dots \oplus V_{k_t})=k_t-k_1+1 $ is called the
length of graduation. A graduation is called of maximum length if $l(\oplus \ll)=dim( \ll)$.
  \end{defn}

  We define the length of an algebra $\ll$ by

 \noindent  $l(\ll)=max \{l(\oplus \ll)  \hbox{ such that } \ll=V_{k_1}\oplus \dots \oplus V_{k_t}\hbox{ is a connected graduation}\}.$

  \begin{defn}\label{def:m.l.}
A Leibniz algebra $\ll$ is called of maximum length if $l(\ll)=dim (\ll)$.
  \end{defn}

\section{Naturally graded quasi-filiform Leibniz algebras}

  The following theorem gives the classification of naturally graded $n$-dimensinal quasi-filiform Lie algebras (see \cite{Liegn}).

  \begin{thm} \label{Liecf} Let $\gg$ be a complex $n$-dimensional naturally
graded quasi-filiform Lie algebra. Then there exists a basis
$\{x_{0}, x_{1}, \dots, x_{n-2}, y\}$ of $\gg$, such that the
multiplication in the algebra has the following form:
$$\small\begin{array}{ll}
L(n, r)\ (n \geq  5,\ 3 \leq  r \leq  2\lfloor\frac{n-1}{2}\rfloor
-1,\ \mbox{r odd}):&Q(n, r) (n \geq  7,\ n \mbox{ odd}, 3 \leq r
\leq  n-4,\ \mbox{r  odd}):
\\[3mm]
\left\{ \begin{array}{l} [x_0 ,x_i ]=x_{i+1} ,\mbox{ }1\le i\le n-3
\\{}
 [x_i ,x_{r-i} ]=(-1)^{i-1}y,\mbox{ }1\le i\le \frac{r-1}{2}
 \end{array} \right.&\left\{ \begin{array}{l}
 [x_0 ,x_i ]=x_{i+1} ,\mbox{ }1\le i\le n-3 \\{}
 [x_i ,x_{r-i} ]=(-1)^{i-1}y,\mbox{ }1\le i\le \frac{r-1}{2} \\{}
 [x_i ,x_{n-2-i} ]=(-1)^{i-1}x_{n-2} ,\mbox{ }1\le i\le \frac{n-3}{2}
 \end{array} \right.
\\[6mm]
\end{array}$$
$$\small\begin{array}{l}
\tau (n, n-3)\ (n \geq 6,\mbox{ n even}):
\\[3mm]
\left\{ {\begin{array}{l} [x_0 ,x_i ]=x_{i+1} ,\mbox{ }1\le i\le n-3
\\{}
[x_{n-1} ,x_1 ]=\frac{(n-4)}{2}x_{n-2} , \\{} [x_i ,x_{n-3-i}
]=(-1)^{i-1}(x_{n-3} +x_{n-1} ),\mbox{ }1\le i\le \frac{n-4}{2} \\{}
[x_i ,x_{n-2-i} ]=(-1)^{i-1}\frac{(n-2-2i)}{2}x_{n-2} ,\mbox{ }1\le
i\le \frac{n-4}{2}
 \end{array}}\right.\\[10mm]
\tau (n, n-4)\ (n \geq
7,\mbox{ n odd}):\\[3mm]
 \left\{ {\begin{array}{l} [x_0 ,x_i ]=x_{i+1}
,\mbox{ }1\le i\le n-3
\\{}
[x_{n-1} ,x_i ]=\frac{(n-5)}{2}x_{n-4+i} ,\mbox{ }1\le i\le 2 \\{}
[x_i ,x_{n-4-i} ]=(-1)^{i-1}(x_{n-4} +x_{n-1} ),\mbox{ }1\le i\le
\frac{n-5}{2} \\{} [x_i ,x_{n-3-i}
]=(-1)^{i-1}\frac{(n-3-2i)}{2}x_{n-3} ,\mbox{ }1\le i\le
\frac{n-5}{2} \\{} [x_i ,x_{n-2-i}
]=(-1)^i(i-1)\frac{(n-3-i)}{2}x_{n-2} ,\mbox{ 2}\le i\le
\frac{n-3}{2}
\end{array}}\right.
\end{array}$$
$$\begin{array}{l}
\small\varepsilon (7, 3):\\
\left\{\begin{array}{ll} [x_0,x_i]=x_{i+1},&1\leq i\leq 3\\{}
[y,x_i]=x_{i+3},&1\leq i\leq 2\\{} [x_1,x_2]=x_{3}+y,&\\{}
[x_1,x_i]=x_{i+1},&3\leq i\leq 4
\end{array}\right.
\end{array}$$
$$\small\begin{array}{ll}
\varepsilon ^{1}(9, 5):& \varepsilon ^{2}(9, 5):\\[3mm]
\left\{\begin{array}{ll} [x_0,x_i]=x_{i+1},&1\leq i\leq 5\\{}
[y,x_i]=2x_{i+5},&1\leq i\leq 2\\{} [x_1,x_4]=x_{5}+y,&\\{}
[x_1,x_5]=2x_{6},& \\{} [x_1,x_6]=3x_{7},& \\{} [x_2,x_3]=-x_{5}-y,&
\\{} [x_2,x_4]=-x_{6},&
\\{} [x_2,x_5]=-x_{7}.&
\end{array}\right.&\left\{\begin{array}{ll} [x_0,x_i]=x_{i+1},&1\leq i\leq 5\\{}
[y,x_i]=2x_{i+5},&1\leq i\leq 2\\{} [x_1,x_4]=x_{5}+y,&\\{}
[x_1,x_5]=2x_{6},& \\{} [x_1,x_6]=x_{7},& \\{} [x_2,x_3]=-x_{5}-y,&
\\{} [x_2,x_4]=-x_{6},& \\{} [x_2,x_5]=x_{7},&\\{} [x_3,x_4]=-2x_{7}.&
\end{array}\right.
\end{array}$$
$$\small\begin{array}{ll}
L_{n-1}\oplus \mathbb{C}\ (n\geq 4):&Q_{n-1}\oplus \mathbb{C}\ (n\geq 7, \ n \ odd):\\[1mm]
\left\{\begin{array}{ll}
[x_0,x_i]=x_{i+1},& 1\leq i\leq n-3.
\end{array}\right.& \left\{\begin{array}{ll}
[x_0,x_i]=x_{i+1},& 1\leq i\leq n-3,\\{}
[x_{i},x_{n-2-i}]=(-1)^{i-1} x_{n-2},& 1\leq i\leq \frac{n-3}{2}.
\end{array}\right.
\end{array}$$
\end{thm}

\

Let $\ll$ be an $n$-dimensional ($n\geq 6$) naturally graded quasi-filiform non Lie Leibniz algebra which has the characteristic
sequence $(n-2,1,1)$ or $(n-2,2)$. The first case (case 2-filiform) has been studied in \cite{Comm.Algebra1} and the second in \cite{J.Lie.Theory2}.
In this work, we will extend this study so we will consider every Leibniz algebras with $C(\ll)=(n-2,2)$.

From the definition of the characteristic sequence it follows the existence of a basis element $e_1 \in \ll\setminus \ll^2$
such that the operator of right multiplication $R_{e_1}$ has one of the following forms:
$$\small\left(
  \begin{array}{c|c}
    J_{n-2} & 0 \\
    \hline 0 &J_{2} \\
  \end{array}
\right)\quad \small\left(
  \begin{array}{c|c}
    J_{2} & 0 \\
    \hline 0 &J_{n-2} \\
  \end{array}
\right)$$

A quasi-filiform Leibniz algebra is called of the type I if the operator $R_{e_1}$ is like the first matrix and of the type II in another case.

 In the following theorems we summarize the results obtained in \cite {NG}.
\begin{thm}\label{Th:NGI}
Let $\mathfrak{A}$ be a naturally graded Leibniz algebra of the first type. Then it is isomorphic to one algebra of the following pairwise non isomorphic families:
$$\small\begin{array}{ll}
\mathfrak{A}^{1, \lambda}:& \mathfrak{A}^{2, \lambda}:\\[2mm]

\begin{cases}
            [y_i,y_1]=y_{i+1}, \quad 1 \leq i \leq n-3\\
            [y_{n-1},y_1]=y_n,\\
            [y_1,y_{n-1}]=\lambda y_{n}.
           \end{cases}&

          \begin{cases}
            [y_i,y_1]=y_{i+1}, \quad 1 \leq i \leq n-3\\
            [y_{n-1},y_1]=y_n,\\
            [y_1,y_{n-1}]=\lambda y_{n}, \quad \lambda \in \{0,1\}\\
            [y_{n-1},y_{n-1}]=y_n.
           \end{cases}
           \\[10mm]
\mathfrak{A}^{3, \lambda}:& \mathfrak{A}^{4, \lambda}:\\[2mm]

\begin{cases}
            [y_i,y_1]=y_{i+1},\quad 1 \leq i \leq n-3\\
            [y_{n-1},y_1]=y_n+y_2,\\
            [y_1,y_{n-1}]=\lambda y_{n}, \quad \lambda \in \{-1,0,1\}.
           \end{cases}&

           \begin{cases}
            [y_i,y_1]=y_{i+1},\quad 1 \leq i \leq n-3\\
            [y_{n-1},y_1]=y_n+y_2,\\
            [y_{n-1},y_{n-1}]=\lambda y_{n}, \quad \lambda \neq 0.
           \end{cases}\\[10mm]
\mathfrak{A}^{5,\lambda, \mu}:&\mathfrak{A}^{6}:\\[2mm]

\begin{cases}
            [y_i,y_1]=y_{i+1},\quad 1 \leq i \leq n-3\\
            [y_{n-1},y_1]=y_n+y_2,\\
            [y_1,y_{n-1}]=\lambda y_{n}, \quad (\lambda, \mu)=(1,1) \hbox{ or } (2,4)\\
            [y_{n-1},y_{n-1}]=\mu y_{n}.
           \end{cases}&

           \begin{cases}
            [y_i,y_1]=y_{i+1},\quad 1 \leq i \leq n-3\\
            [y_{n-1},y_1]=y_n,\\
            [y_1,y_{n-1}]=-y_{n}, \\
            [y_{n-1},y_{n-1}]= y_{2},\\
            [y_{n-1},y_{n}]= y_{3}.
           \end{cases}
\end{array}$$
\end{thm}

\begin{thm}\label{Th:NGII}
Let $\mathfrak{A}$ be a naturally graded Leibniz algebra of the second type. Then it is
isomorphic to one algebra of the following pairwise non isomorphic families:

$n$ even

$$\small\begin{array}{ll}
\mathfrak{A}^{1}:&\mathfrak{A}^{2}:\\[2mm]

\begin{cases}
            [y_1,y_1]=y_{2},\\
            [y_{i},y_1]=y_{i+1}, \quad 3 \leq i \leq n-1\\
            [y_1,y_{i}]=-y_{i+1}, \quad 3 \leq i \leq n-1.
           \end{cases}& \begin{cases}
            [y_1,y_1]=y_{2},\\
            [y_{i},y_1]=y_{i+1}, \quad 3 \leq i \leq n-1\\
            [y_1,y_3]=y_2-y_4,\\
            [y_1,y_{j}]=-y_{j+1}, \quad 4 \leq j \leq n-1.
           \end{cases}
        \end{array}$$

$$\small\begin{array}{ll}
\mathfrak{A}^{3}:&\mathfrak{A}^{4}:\\[2mm]

\begin{cases}
            [y_1,y_1]=y_{2},\\
            [y_{i},y_1]=y_{i+1}, \quad 3 \leq i \leq n-1\\
            [y_3,y_3]=y_2,\\
            [y_1,y_{i}]=-y_{i+1}, \quad 3 \leq i \leq n-1.
           \end{cases}& \begin{cases}
            [y_1,y_1]=y_{2},\\
            [y_{i},y_1]=y_{i+1}, \quad 3 \leq i \leq n-1\\
            [y_1,y_3]=2y_2-y_4,\\
            [y_3,y_3]=y_2,\\
            [y_1,y_{j}]=-y_{j+1}, \quad 4 \leq j \leq n-1.
           \end{cases}
\end{array}$$

$n$ odd, $\mathfrak{A}^{1},\mathfrak{A}^2,\mathfrak{A}^3, \mathfrak{A}^4$

$$\small\begin{array}{ll}
\mathfrak{A}^5:&\mathfrak{A}^{6, \lambda}:\\[2mm]
\begin{cases}
            [y_1,y_1]=y_{2},\\
            [y_{i},y_1]=y_{i+1}, \quad 3 \leq i \leq n-1\\
            [y_1,y_{i}]=-y_{i+1}, \quad 4 \leq i \leq n-1\\
            [y_i,y_{n+2-i}]=(-1)^iy_n, \quad 3 \leq i \leq n-1.
           \end{cases}& \begin{cases}
            [y_1,y_1]=y_{2},\\
            [y_{i},y_1]=y_{i+1}, \quad 3 \leq i \leq n-1\\
            [y_1,y_3]=\lambda y_2-y_4, \quad \lambda \in \{1,2\}\\
            [y_1,y_{j}]=-y_{j+1}, \quad 4 \leq j \leq n-1\\
            [y_i,y_{n+2-i}]=(-1)^iy_n, \quad 3 \leq i \leq n-1.
           \end{cases}\\[19mm]
\mathfrak{A}^{7,\lambda}:&\mathfrak{A}^{8, \lambda, \mu}:\\[2mm]

\begin{cases}
            [y_1,y_1]=y_{2},\\
            [y_{i},y_1]=y_{i+1}, \quad 3 \leq i \leq n-1\\
            [y_3,y_3]=\lambda y_2, \quad \lambda \neq 0\\
            [y_1,y_{i}]=-y_{i+1}, \quad 3 \leq i \leq n-1\\
            [y_i,y_{n+2-i}]=(-1)^iy_n, \quad 3 \leq i \leq n-1.
           \end{cases}&\begin{cases}
            [y_1,y_1]=y_{2},\\
            [y_{i},y_1]=y_{i+1}, \quad 3 \leq i \leq n-1\\
            [y_1,y_3]=\lambda y_2-y_4, \\
            [y_3,y_3]=\mu y_2, \\
            [y_1,y_{j}]=-y_{j+1}, \quad 4 \leq j \leq n-1\\
            [y_i,y_{n+2-i}]=(-1)^iy_n, \quad 3 \leq i \leq n-1.
           \end{cases}\\[3mm]
&\mbox{with } (\lambda, \mu)=(-2,1),(2,1) \hbox{ or }(4,2)
\end{array}$$
\end{thm}

The study of Leibniz algebras of maximum length can be simplified by following the reasoning in the
proof of theorems \ref{Th:NGI} and \ref{Th:NGII}, see \cite{NG}. The next proposition gives
the structure of a naturally graded $n$-dimensional Leibniz algebra:
\begin{pr} \label{pr:NG'} Let $\ll $ be a naturally graded quasi-filiform Leibniz algebra, then it is
isomorphic to one algebra of the pairwise non isomorphic families:\\

$\begin{array}{ll}
NG1: & NG2: \\[2mm]
\begin{cases}
                [e_1,e_1]=e_2,  \\
                [e_{i},e_{1}]=e_{i+1},  &3\leq i \leq n-1 \\
                [e_{1},e_{3}]=\lambda e_2-e_4, \\
                [e_3,e_3]=\mu e_2,& \lambda, \ \mu \in \mathbb{C}\\
                [e_1,e_i]=-e_{i+1}, &4\leq i \leq n-1.
          \end{cases} & \begin{cases}
                [e_1,e_1]=e_2, \\
                [e_{i},e_{1}]=e_{i+1},  &3\leq i \leq n-1 \\
                [e_{1},e_{3}]=\lambda e_2-e_4, \\
                [e_3,e_3]=\mu e_2,& \lambda, \ \mu \in \mathbb{C}\\
                [e_1,e_i]=-e_{i+1},  &4\leq i \leq n-1\\
                [e_i,e_{n+2-i}]=(-1)^ie_n,  &3 \leq i \leq n-1.\\
          \end{cases} \\[19mm]
 \end{array}$\\

$ \begin{array}{ll}
NG3: & NG4: \\[2mm]
\begin{cases}
                [e_i,e_1]=e_{i+1}, &1\leq i \leq n-3 \\
                [e_{n-1},e_{1}]=e_{n}+\alpha e_2, \\
                [e_{1},e_{n-1}]=\beta e_n, \\
                [e_{n-1},e_{n-1}]=\gamma e_n, & \alpha, \ \beta, \ \gamma \in \mathbb{C}.
          \end{cases} & \begin{cases}
                [e_i,e_1]=e_{i+1}, & 1\leq i \leq n-3 \\
                [e_{n-1},e_{1}]=e_{n},  \\
                [e_{1},e_{n-1}]=-e_n, \\
                [e_{n-1},e_{n-1}]=e_2,\\
                [e_{n-1},e_n]=e_3.
          \end{cases}
   \end{array}$

\noindent where NG1 and NG2 correspond to algebras of the type II and NG3 and NG4 to algebras of the type I.
\end{pr}

\section{Quasi-filiform Leibniz algebras of maximum lenght}

Let $\ll$ be an $n$-dimensional ($n\geq 6$) quasi-filiform non Lie Leibniz algebra with the characteristic
sequence $(n-2,2)$. The proposition \ref{pr:NG} shows a structure of these algebras and requires the previous result.
\begin{pr}\label{pr:NG}
Let $\ll$ be an $n$-dimensional quasi-filiform Leibniz algebra. Then $\ll$ is isomorphic to one algebra of the following pairwise non isomorphic families:

$$\small\begin{array}{l}
\widetilde{NG1}:\\[3mm]
\begin{cases}
                [e_1,e_1]=e_2,\\
                [e_{i},e_{1}]=e_{i+1},  \quad 3\leq i \leq n-1 \\
                [e_{1},e_3]=\lambda e_2-e_4+ \mu_n e_{n},  \\
                [e_3,e_3]=\mu e_2+\gamma_n e_n, \\
                [e_1,e_i]=-e_{i+1}+\gamma_{i,n}e_n,  \quad 4\leq i \leq n-1\\
                [e_i,e_j]=(*)e_{i+j-1}+(*)e_{i+j}+\dots+(*)e_n, \quad \forall (i,j)\neq (1,1),(1,3),(3,3),(3,1),(i,1),(1,i)\\
                [e_2,e_i]=z_ie_n, \quad i \neq 2,n.
          \end{cases}\\[15mm]

\\

\widetilde{NG2}:\\[3mm]
\begin{cases}
                [e_1,e_1]=e_2+\alpha_ne_n,\\
                [e_{i},e_{1}]=e_{i+1},  \quad 3\leq i \leq n-1 \\
                [e_{1},e_3]=\lambda e_2-e_4+(*)e_{n},  \\
                [e_3,e_3]=\mu e_2+ (*)e_n, \\
                [e_1,e_i]=-e_{i+1}+(*)e_n,  \quad 4\leq i \leq n-1\\
                [e_i,e_{n+2-i}]=(-1)^ne_n,  \quad 3 \leq i \leq n-1\\
                [e_i,e_j]=(*)e_{i+j-1}+ \dots+(*)e_n, \quad \forall (i,j)\neq (1,1),(1,3),(3,3),(3,1),(k,n+2-k)\\
                \qquad \qquad\qquad \quad  \qquad \qquad \qquad \qquad 3 \leq k\leq n-1.
          \end{cases}\\[15mm]

\\


\widetilde{NG3}:\\[3mm]
\begin{cases}
                [e_i,e_1]=e_{i+1}+(*)e_{i+2}+ \dots+ (*)e_{n-2}, &1\leq i \leq n-3 \\
                [e_{n-1},e_{1}]=e_{n}+\alpha e_2+(*)e_3+ \dots+(*)e_{n-2}, \\
                [e_{1},e_{n-1}]=\beta e_n+(*)e_3+ \dots+(*)e_{n-2}, \\
                [e_{n-1},e_{n-1}]=\gamma e_n+(*)e_3+ \dots+(*)e_{n-2},\\
                [e_i,e_{n-1}]=(*)e_{i+2}+\dots + (*)e_{n-2}, &2 \leq i \leq n-3\\
                [e_n, e_{n-1}]=(*)e_4 + \dots+(*)e_{n-2}.
          \end{cases}
                    \end{array}$$
$$\small\begin{array}{l}

\widetilde{NG4}:\\[4mm]
\begin{cases}
          [e_i,e_1]=e_{i+1}+(*)e_{i+2}+ \dots+(*)e_{n-2}, &1\leq i \leq n-3 \\
          [e_{n-1},e_{1}]=e_{n}+(*)e_3+\dots+(*)e_{n-2},  \\
          [e_{1},e_{n-1}]=-e_n+(*)e_3+\dots+(*)e_{n-2}, \\
          [e_{n-1},e_{n-1}]=e_2+(*)e_3+\dots+(*)e_{n-2},\\
          [e_{n-1},e_n]=e_3+(*)e_4+\dots+(*)e_{n-2},\\
          [e_i,e_j]=(*)e_{i+j+1}+ \dots+(*)e_n, &2\leq i,j\leq n-2\\
          [e_i,e_{n-1}]=(*)e_{i+2}+\dots+(*)e_{n-2}, & 2 \leq i \leq n-2\\
          [e_n,e_{n-1}]=(*)e_4+\dots+(*)e_{n-2},\\
          [e_j,e_n]=(*)e_5 +\dots +(*)e_{n-2}, & j=2,5\\
          [e_i,e_n]=(*)e_{i+3}+ \dots+(*)e_{n-2}, & i=1 \wedge 3\leq i \leq n-2\\
          [e_n,e_n]=(*)e_5+\dots+(*)e_{n-2}.
          \end{cases}
          \end{array}$$
          where the asterisks (*) denote appropiate coefficients at the basis elements.
\end{pr}

\begin{dem}
These families are obtained from the structure of the naturally graded quasi-filiform Leibniz algebras (see Proposition \ref{pr:NG'}) considering its natural graduations.
\end{dem}

Since the study on quasi-filiform Lie algebras of maximum length has been done in \cite{Reyes1}, we shall
 restrict ourselves to Leibniz algebras which are non Lie. Moreover, the classification of 2-filiform Leibniz
 algebras of maximum length is given in \cite{J.Lie.Theory2}. Thus, the following results close the study of maximum
 length for Leibniz algebras with nilindex $n-2$, where $n$ is the dimension of $\ll$ and verified $n \geq 6$.

\begin{thm}\label{th:tipo I}
Let $\ll$ be an $n$-dimensional quasi-filiform non Lie Leibniz algebra of maximum length of
the type I. Then the algebra $\ll$ is isomorphic to one algebra of the following pairwise non isomorphic families:
$$\begin{array}{ll}
M^{1,\delta}:& M^{2,\lambda}:\\[3mm]
\begin{cases}
            [y_1,y_1]=y_n,\\
            [y_{n-1},y_1]=y_2,\\
            [y_i,y_1]=y_{i+1}, \quad 2 \leq i \leq n-3\\
            [y_{n-1},y_{n-1}]=\delta y_4, \quad \delta \in \{0,1\}\\
            [y_i,y_{n-1}]=\delta y_{3+i},\quad 2\leq i\leq n-5.\\
           \end{cases} &
                       \begin{cases}
                 [y_i,y_1]=y_{i+1}, \quad 1\leq i \leq n-3\\
                 [y_{n-1},y_1]= y_n,\\
                 [y_1,y_{n-1}]=\lambda y_{n}, \quad \lambda \in \mathbb{C}.
                \end{cases}
                \end{array}$$
\end{thm}
\begin{dem}
Let $<e_1,e_2, \dots, e_n> $ be the basis used in proposition \ref{pr:NG}. We can consider the following cases:

 \noindent \underline{Case 1:} Let $e_n \notin R(\ll)$ be.

 In this cases we have the family $\widetilde{NG4}$, where $\{e_2, e_3, \dots,e_{n-2}\}\subseteq R(\widetilde{NG4})$.

  In order to study the length of the algebras from the family $\widetilde{NG4}$, let us take
$$x_{s}=e_{1}+ \sum_{i=2}^{n} a_{i}e_{i} \quad  \quad x_t=e_{n-1}+\sum_{k=1; k \neq n-1}^{n}b_k e_k, \hbox{ where } a_{n-1}b_1 \neq 1$$
and
\begin{align*}
[x_s,x_s]&=(1+a_{n-1}^2)e_2+(*)e_3+\dots+(*)e_n, \\
[x_t,x_t]&=(1+b_1^2)e_2+(*)e_3+ \dots+(*)e_n,
\end{align*}
\begin{align*}
[x_t,x_s]&=(b_1+a_{n-1})e_2+(*)e_3+ \dots +(*)e_{n-1}+(1-a_{n-1}b_1)e_n,\\
[x_s,x_t]&=(b_1+a_{n-1})e_2+(*)e_3+ \dots +(*)e_{n-1}+(a_{n-1}b_1-1)e_n.
\end{align*}

We will assume two subcases:

\noindent {\bf Subcase 1.1: }If $1+a_{n-1}^2 \neq 0$. Let us consider
$$\underbrace{[[x_s,x_s],x_s],...,x_s}_{i\hbox{ } times}]=(1+a_{n-1}^2)e_i+(*)e_{i+1}+\dots+(*)e_n$$
and the new basis $\{y_1,\dots, y_n\}$, defined as:
$$y_1=x_s,\quad y_i=[y_{i-1},x_s], \ \hbox{ with } \ 2\leq i \leq n-2, \quad y_{n-1}=x_t,\quad y_{n}=[x_s,x_t]$$
with the graduation of maximum length: $V_{k_{s}} \oplus V_{2k_{s}}\oplus \dots\oplus V_{(n-2)k_s}\oplus V_{k_t}\oplus V_{k_s+k_t}$.
This basis verifies the following product:
$$\begin{array}{ll}
[y_i,y_1]=y_{i+1}, & 1 \leq i \leq n-3\\{}
[y_1,y_{n-1}]=y_{n}.&
\end{array}$$

Furthermore, $$ y_n \notin  R(\widetilde{NG4}),\ [y_1,y_{n-1}]+[y_{n-1},y_1] \in R(\widetilde{NG4}), \mbox{ and } \ [y_{n-1},y_1]=A y_n$$
then $A=-1$. Let us check the others products.

We shall consider
$$\begin{cases}
[y_{n-1},y_1]=[x_t,x_s]=(b_1+a_{n-1})e_2+(*)e_3+ \dots +(*)e_{n-1}+(1-b_1a_{n-1})e_n,\\
[y_1,y_{n-1}]=[x_s,x_t]=(b_1+a_{n-1})e_2+(*)e_3+ \dots +(*)e_{n-1}+(b_1a_{n-1}-1)e_n,\\
[y_{n-1},y_1]=-[y_1,y_{n-1}],
 \end{cases}$$
therefore $b_1=-a_{n-1}$.

 From $a_{n-1}=-b_1$  and  $a_{n-1}b_1 \neq 1$, one has  $a_{n-1}^2 \neq -1$, and so $[y_{n-1},y_n]=$ $=-(a_{n-1}^2+1)e_3+(*)e_4+\dots+(*)e_n$ is obtained. Hence $[y_{n-1},y_n]=B y_3$ where $B \neq 0$.
Therefore
$$\begin{cases}
[y_{n-1},y_n]=B y_3, & B\neq 0,\\
[y_{n-1},y_{n}] \subseteq V_{2k_t+k_s},\\
y_3 \subseteq V_{3k_s}.
\end{cases}$$

 Thus, we prove that $k_s=k_t$, which is a contradiction of the definition of maximum length.

\

\noindent {\bf Subcase 1.2: }If $1+ a_{n-1}^2=0 $. Since
$$\begin{cases}
[x_t,x_s]=(b_1+a_{n-1})e_2+(*)e_3+ \dots +(*)e_{n-1}+(1-b_1a_{n-1})e_n=A y_n \neq 0,\\
[x_s,x_t]=(b_1+a_{n-1})e_2+(*)e_3+ \dots +(*)e_{n-1}+(b_1a_{n-1}-1)e_n= y_n \neq 0,\\
[x_t,x_s]+[x_s,x_t] \in R(\widetilde{NG4}),
\end{cases}$$
and by the property of maximum length we have $(A+1)y_n=0$, i.e, $A=-1$, hence $b_1=-a_{n-1}$ and $a_{n-1} b_1\neq 1,$ it implies that $a^2_{n-1}\neq -1$.
 But it is a contradiction of the hypothesis $1+ a_{n-1}^2=0 $.

\

Thus, we can conclude the algebra $\widetilde{NG4}$ does not have maximum length.

\

\noindent \underline{Case 2:} Let $e_n \in R(\ll)$ be.

In this case we have the family $\widetilde{NG3}$, where $\{e_2,e_3, \dots,e_{n-2},e_n\} \subseteq  R(\widetilde{NG3})$.

Let us take the previous basis, we have the products:
$$\small (1)\left\{\begin{array}{l}
[x_s,x_s]=(1+a_{n-1}\alpha)e_2+(*)e_3+\dots+(*)e_{n-2}+(a_{n-1}\beta +a_{n-1}+a_{n-1}^2\gamma )e_n, \\{}
[x_t,x_t]=(b_1^2+b_1\alpha)e_2+(*)e_3+ \dots+(*)e_{n-2}+(b_1+b_1\beta + \gamma)e_n,\\{}
[x_t,x_s]=(b_1+\alpha)e_2+(*)e_3+ \dots +(*)e_{n-2}+(1+b_1a_{n-1}\beta+a_{n-1}\gamma)e_n,\\{}
[x_s,x_t]=(b_1+a_{n-1}b_1\alpha)e_2+(*)e_3+ \dots +(*)e_{n-2}+(b_1a_{n-1}+\beta+a_{n-1}\gamma)e_n,\\{}
\underbrace{[[x_s,x_s],x_s]...,x_s}_{i-times}]=(1+a_{n-1}\alpha)e_i+(*)e_{i+1}+\dots+(*)e_{n-2}.
\end{array}\right.$$

\noindent {\bf Subcase 2.1: }If $1+a_{n-1}\alpha=0$.\\

If $b_1(b_1+\alpha) \neq 0$, then in a similar way to previous cases, we can prove that there is not any algebra of maximum length.
Therefore we shall assume that $b_1(b_1+\alpha)=0$.

           Note that indeed if $b_1+\alpha=0,$ then we will have $1+a_{n-1}\alpha=0$ and hence $a_{n-1}b_1=1$, which would
            be impossible. Therefore $b_1=0$.
          In this case we obtain
          \begin{align*}
         [x_s,x_s]&=(*)e_3+\dots+(*)e_{n-2}+(a_{n-1}\beta +a_{n-1}+a_{n-1}^2\gamma )e_n, \\
         [x_t,x_t]&=(*)e_3+ \dots+(*)e_{n-2}+\gamma e_n,\\
         [x_t,x_s]&=\alpha e_2+(*)e_3+ \dots +(*)e_{n-2}+(1+a_{n-1}\gamma)e_n,\\
         [x_s,x_t]&=(*)e_3+ \dots +(*)e_{n-2}+(\beta+a_{n-1}\gamma)e_n,\\
         [[[x_t,\underbrace{x_s],x_s],...,x_s]}_{(i-1)-\hbox{times}}&=\alpha e_i+(*) e_{i+1}+\dots+(*) e_{n-2},\quad 3\leq i\leq n-2.
       \end{align*}

       Since $[x_s,x_t]=D[x_t,x_s]$, it implies that $D\alpha=0$ as $\alpha\neq 0$  we have $D=0$.
       Thus, $\beta+a_{n-1}\gamma=0$ is achieved and $[y_1,y_{n-1}]=0$.

       Put $y_1=x_s,$ $y_2=[x_t,x_s],$ $ y_{n-1}=x_t,$ $ y_i=[y_{i-1},y_1]$ with $3\leq i \leq n-2 $ and $ y_n=[y_1,y_1].$

       This basis has the following graduation of maximum length: $V_{k_{s}} \oplus V_{k_t+k_{s}}\oplus V_{k_t+2k_s}\oplus\dots\oplus V_{k_t+(n-3)k_s}\oplus V_{k_t}\oplus V_{2k_s}.$

       We have that $[y_{n-1},y_{n-1}]=[x_t,x_t]=(*)e_3+ \dots+(*)e_{n-2}+\gamma e_n=Ay_m$, where $m \notin \{1,2,n-1\}.$
               On the other hand $[y_{n-1},y_{n-1}] \in V_{2k_t}.$ Therefore $2k_t \in \{2k_s,k_t+(m-1)k_s\}$ with $3 \leq m \leq n-2.$ If $2k_t=2k_s$ then this is a contradiction of the definition of maximum length. As $y_m=\alpha e_m+(*)e_{m+1}+ \dots + (*)e_{n-2}$ with $3 \leq m \leq n-2$ then $\gamma =0$. Hence the form $\beta + a_{n-1}\gamma =0$, we get $\beta =0.$
               
               The cases $m=3$ and $5 \leq m \leq n-3$ are impossible considering
                the connectedness of the graduation of maximum length. Thus let us see $m=4$ and $m=n-2.$
               
               \
               
               Let $m=n-2,$ and $n >6$. It is impossible to find a connected graduation of maximum length. If $n=6$, then we have $m=4$.

\

               Let $m=4.$ We have $[y_{n-1},y_{n-1}]=\delta y_4$ and by using the Leibniz identity $[y_i,y_{n-1}]=\delta y_{i+1}$ with $2 \leq i \leq n-5$ and $[y_n,y_{n-1}]=0$ are achieved.

      In order to find the table of multiplications of the algebra it is enough to show $[y_n,y_1]$ since $\{y_2,y_3, \dots,y_{n-2},y_n\} \subseteq R(\ll)$ and $y_{n-2} \in Cent(\ll).$ By the construction, $[y_n,y_1]=(1+a_{n-1}\alpha )e_3+(*)e_4+ \dots +(*)e_{n-2}$. As $1+a_{n-1}\alpha \neq 0$ it implies that $[y_n,y_1]=By_3$.
      Hence, by the properties of the graduation, $k_t=k_s$ which is a contradiction of the definition of maximum length. Thus $[y_n,y_1]=0.$

      Note that if $\delta \neq 0$, we can consider $\delta=1$ using an easy change of basis. Therefore, the algebra $M^{1,\delta}$ with $\delta=0$ or $\delta=1$
       is achieved. $M^{1,0}$ and $M^{1,1}$ are non isomorphic because $dim(R(M^{1,0}))\neq dim(R(M^{1,1})).$

        To see its  maximum length it is enough to consider the graduation
       $V_1\oplus V_2 \oplus \dots \oplus$ $ \oplus V_n$, where $V_1=<y_1>,\ $ $V_2=<y_n>,\ $ $V_3=<y_{n-1}>$
       and $V_i=<y_{i-2}> $  with $ 4 \leq i \leq n.$

\

\noindent {\bf Subcase 2.2: } If $1+a_{n-1}\alpha \neq 0$.\\

Using the multiplications of $(1)$ we can distinguish the following cases:

\

1) If $b_1+b_1\beta+\gamma \neq 0\ \wedge \ det\small\left(
                     \begin{array}{cc}
                     1+a_{n-1}\alpha & a_{n-1}\beta+a_{n-1}+a_{n-1}^2\gamma\\
                     b_1^2+b_1\alpha & b_1+b_1\beta+\gamma\\
                     \end{array}
                    \right) \neq 0$\\

Since the determinant is not equal to zero, we have  that
$[x_s,x_s]$ and $[x_t,x_t]$ are linearly independent. We define a
new basis: $y_1=x_s,$  $y_i=[y_{i-1},y_1]$ with $2\leq i\leq
n-2$, $y_{n-1}=x_t$ and $y_n=[x_t,x_t]$. Making the following
graduation: $V_{k_s}\oplus V_{2k_s}\oplus \dots \oplus
V_{(n-2)k_s}\oplus V_{k_t}\oplus V_{2k_t}$ and if $[x_t,x_s]\neq
0$ then we have that $[x_t,x_s]\in V_{k_t+k_s}$ and
$k_t+k_s\notin \{k_s,k_t,2k_s,2k_t\}$. It implies $[x_t,x_s]= A
y_m$ with $A\neq 0$ and $3\leq m\leq n-2$ $\Rightarrow$
$k_t=(m-1)k_s$. If $k_s>0$ then $2k_s\leq k_t\leq (n-2)k_s$ and
that is a contradiction. Thus, $[x_t,x_s]=0$ hence $b_1+\alpha=0$
and $1+a_{n-1}b_1\beta+a_{n-1}\gamma=0$. Analogously for
$[x_s,x_t]$, we have $b_1+a_{n-1}b_1\alpha=0$ and
$a_{n-1}b_1+\beta +a_{n-1}\gamma=0$.

We have $b_1=\alpha=0,$ $\beta=1,$ $a_{n-1}\gamma\neq 0$. We know that $\{y_2,\dots, y_{n-2},y_n\}\subseteq R(\widetilde{NG3})$ and
the products $[y_i,y_1]=y_{i+1}$ with $1\leq i\leq n-3$ and $[y_{n-1},y_{n-1}]=y_n$.
Using the properties of the graduation $[y_{n-1},y_1],\ [y_{1},y_{n-1}]$ and $[y_{n},y_{n-1}]$ are obtained and it is archived that there is not any Leibniz algebras of maximum length.

\

2) If $a_{n-1}b_1+\beta+a_{n-1}\gamma \neq 0 \ \wedge \ det\small\left(
                     \begin{array}{cc}
                     1+a_{n-1}\alpha & a_{n-1}\beta+a_{n-1}+a_{n-1}^2\gamma\\
                     b_1(1+a_{n-1}\alpha) & b_1a_{n-1}+\beta+a_{n-1}\gamma\\
                     \end{array}
                    \right) \neq 0,$\\
     then we can take the basis $\{y_1, \dots, y_n\}$, defined as:
  $y_1=x_s,$ $y_i=[y_{i-1},y_1],$ with $2\leq i \leq n-2,$ $y_{n-1}=x_t,$ and $y_n=[y_1,y_{n-1}],$
          and its graduation of maximum length: $V_{k_{s}}\oplus V_{2k_s} \oplus \dots \oplus V_{(n-2)k_s}\oplus V_{k_{t}}\oplus V_{k_{t}+k_s}.$\\

          \noindent Since the determinant does not vanish, we get $\beta+a_{n-1}\gamma \neq 0$.
          Therefore, $\{e_2,$ $\dots,$ $e_{n-2},$ $e_n\} \subseteq R(\ll)$ and $e_{n-2}\in Cent(\ll)$,
           hence $\{y_2, \dots, y_{n-2},y_n\} \subseteq R(\ll)$, and $y_{n-2}\in Cent(\ll)$. Let us define $R_{y_1}$ and $R_{y_{n-1}}$.

\

           \noindent If $[y_n,y_1]\neq 0$ then $[y_n,y_1]\in V_{2k_s+k_t}$. Therefore $2k_s+k_t\notin \{2k_s,k_t,k_t+k_s, 3k_s\}$ because
            $k_s \neq 0 \neq k_t$ and $k_s \neq k_t$. Moreover, it is evident that $2k_s+k_t \neq k_s$, otherwise $k_t=-k_s$ and this implies
             $[y_n,y_1]=Ay_1 $ with $ A \neq 0$, but it is a contradiction since $[y_n,y_1]\in L_{3}$ and $y_1 \in L_1$. On the other hand we can assume
            $\ 2k_s+k_t \neq mk_s$ for $4\leq m \leq n-2,$ (otherwise we will obtain a contradiction of maximum length).
            As a result we can conclude that $[y_n,y_1]=0$.

\

\noindent A similar study for $[y_i,y_{n-1}],$ with $2 \leq i \leq n-3$ permits to
prove that $[y_i,y_{n-1}]=0 \quad 2 \leq i \leq n-3$, $[y_n,y_{n-1}]=0,$ $[y_{n-1},y_{n-1}]=0 $ and so $b_1(b_1+\alpha)=b_1+b_1\beta+\gamma=0$.
Moreover, we have
     $$\begin{cases}
         [y_1,y_{n-1}]=B [y_{n-1},y_1],\\
         [y_1,y_{n-1}]=b_1(b_1+a_{n-1}b_1\alpha)e_2+(*)e_3+ \dots + (*)e_n, \\
         [y_{n-1},y_1]=(b_1+\alpha)e_2+(*)e_3+ \dots+ (*)e_n.
         \end{cases}$$
         it implies that $b_1=0, \gamma=0 \hbox{ and } B=\beta \neq 0$. Thus $\alpha=0$.\\

         Finally we have found the family defined as:
         \begin{align*}
         [y_i,y_1]&=y_{i+1}, & 1 \leq i \leq n-3\\
         [y_{n-1},y_1]&=B y_n, & B \in \mathbb{C}\setminus \{0\}\\
         [y_1,y_{n-1}]&=y_n.
         \end{align*}

         Taking $y^{'}_n=B y_n$, the family $M^{2,\lambda}$ is achieved, with $\lambda \in \mathbb{C}\setminus \{0\}$.
         To get maximum length it is enough to consider the following graduation:
         $$\ll=V_1 \oplus V_2 \oplus \dots \oplus V_n, \quad \hbox{ where }\quad  V_i=<y_i> \quad 1 \leq i \leq n.$$

\

3)    If $1+a_{n-1}b_1\beta+a_{n-1}\gamma \neq 0 \  \wedge \ det\small\left(
                     \begin{array}{cc}
                     1+a_{n-1}\alpha & a_{n-1}\beta+a_{n-1}+a_{n-1}^2\gamma\\
                     b_1+\alpha & 1+ b_1a_{n-1}\beta+a_{n-1}\gamma\\
                     \end{array}
                    \right) \neq 0$

   By following similar arguments as in the previous cases, we obtain the maximum length family $M^{2,\lambda}, $ with $\lambda\in \mathbb{C}$.

   Note that $M^{1,0}$ and $M^{2,\lambda}$ with $\lambda \neq 0$ are non isomorphic because $dim(R(M^{1,0})) \neq dim(R(M^{2,\lambda}))$. Analogously, $M^{1,0} \ncong M^{2,0}$ can be proved using change of basis and $M^{1,1} \ncong M^{2,\lambda}$ with $\lambda \in \mathbb{C}$ because $dim(\mathcal{I}(M^{1,1}))
    \neq dim(\mathcal{I}(M^{2,\lambda}))$.
\end{dem}

\begin{thm}\label{th:tipo II}
Let $\ll$ be an $n$-dimensional quasi-filiform non Lie Leibniz algebra of maximum length of the type II. Then the algebra $\ll$ is isomorphic to
one algebra of the family

$M^{3,\alpha}:\begin{cases}
                   [y_1,y_1]=y_2,\\
                   [y_i,y_1]=y_{i+1}, &3\leq i \leq n-1,\\
                   [y_1,y_i]=-y_{i+1}, &4\leq i\leq n-1,\\
                   [y_3,y_3]=\alpha y_6, & \alpha=0 \hbox{ if } n> 6,\  \alpha \in \{0,1\} \hbox{ if } n=6.
                 \end{cases} $
\end{thm}

\begin{dem}
Let $\{e_1,e_2, \dots, e_n\} $ be the basis used in the proposition \ref{pr:NG}. In this case we have $e_4 \notin R(\ll)$.

 Let us consider two cases:
 \begin{itemize}
  \item {Case 1: } Let $n$ be even.

  In this case we have the algebra $\widetilde{NG1}$ and we are going to study its length. Analogously to the previous theorem, if $n\neq 6$, then we obtain $M^{3,0}$ (which has maximum length)  and if $n=6$ the family of maximum length $M^{3,\alpha}$ with $\alpha\in \mathbb{C}$.

  \item {Case 2: } Let $n$ be odd.

  In the same way, when we study the length of $\widetilde{NG2}$, we can prove that any maximum length algebra does not exist  in this family.
 \end{itemize}
\end{dem}

It remains to prove that if we consider the extension of naturally graded quasi-filiform Lie algebras, then we obtain
Lie algebras too (see Theorem \ref{Lieext}) and their study of maximum length can be found in \cite{Reyes1}.
The next theorem and their corollary will be used to prove the Theorem \ref{Lieext}.

\begin{thm}\label{Lie}
 Let $\ll$ be an $n$-dimensional Leibniz algebra and let $\{y_0,y_1,\dots,y_{n-1}\}$ be a basis. Let $\{y_0,y_1, \dots,y_s\}$ be
 the generators of $\ll$. If the following property is true
$ [y_i,y_j]=-[y_j,y_i] $ for all $y_j \in \ll$, with $0\leq i\leq s$, then $\ll$ is a Lie algebra.
\end{thm}

\begin{dem}
It is enough to prove that the rest of the multiplication is antisymmetric, i.e., to prove
that $[y_i,y_m]=-[y_m,y_i] $ where $ i,m \in \{s+1,\dots,n-1\}$. Let us consider two steps
\begin{itemize}
\item{ Step 1:} If $y_m=[y_{i_0},y_{j_0}] \hbox{, where } i_0,j_0 \in \{0,1, \dots, s\}$, then
from the Leibniz identity and the hypothesis of the theorem we have:
$[y_i,y_m]=-[y_m,y_i]$ with  $0 \leq i \leq n-1$.

\item{ Step 2: } If $y_m=[y_{p},y_{i_0}] \hbox{ with } i_0 \in \{0,1, \dots, s\}$ and  $p \in \{s+1, \dots, n-1\}$, then
from the Leibniz identity, the hypothesis of the theorem and the step before we have:
$[y_i,y_m]=-[y_m,y_i]$ with $0\leq i\leq n-1$.
\end{itemize}
Note that $y_i$ is an arbitrary element of $ \ll$.
\end{dem}

In order to make the proof of Theorem \ref{Lieext} easier, we will need the following corollaries.
\begin{cor}\label{Lema1}
 Let $\ll$ be an $n$-dimensional Leibniz algebra  and $\mathcal{B}=\{y_0,$ $y_1,$ $\dots,$ $y_{n-1}\}$ a basis. Let $\{y_0,y_1, \dots,y_s\}$ be
 the generators of $\ll$. If $\mathcal{B}$ satisfies the following conditions:
 \begin{align*}
 i) \quad &[y_i,y_j]=-[y_j,y_i] \quad \forall i,j \in \{0,1, \dots,s\},\\
ii) \quad &y_i=[y_{i_0},y_{i-1}] \quad s+1 \leq i \leq n-1,\\
 iii) \quad &[y_i,y_{i_0}]=-[y_{i_0},y_i] \quad 1\leq i_0 \leq s, \quad s+1\leq i\leq n-1.
\end{align*}
then $\ll $ is a Lie algebra.

\end{cor}

\begin{dem}
We shall prove $[y_i,y_m]=-[y_m,y_i] $ with $s+1\leq i\leq n-1$, which is directly obtained from the above theorem.
\end{dem}

\begin{cor}\label{Lema2}
 Let $\ll$ be an $n$-dimensional  Leibniz algebra and $\{y_0,y_1,\dots,y_{n-1}\}$ a basis. Let $\{y_{0},y_{1}\}$ be the generators of $\ll$. If the following assertions are achieved
 \begin{align*}
   i)\quad &[y_i,y_j]=-[y_j,y_i] \quad \forall i,j \in \{0,1\},\\
ii)\quad & \begin{cases}
         y_i=[y_{0},y_{i-1}] \quad 2 \leq i \leq n-1,\\
         y_{n-1}=[y_{1},y_{p}]  \quad 2\leq p \leq n-2.
     \end{cases}\\
 iii) \quad &\begin{cases}
      [y_i,y_0]=-[y_0,y_i]  \quad 2\leq i\leq n-1\\
      [y_{n-1},y_1]=-[y_1,y_{n-1}].
     \end{cases}
\end{align*}
then $\ll $ is a Lie algebra.
\end{cor}

\begin{dem}
It is sufficient to show the following products
\begin{itemize}
\item{ Let us see  $[y_{1},y_i]=-[y_i,y_1] $ with $ 2\leq i \leq n-2.$}\\
Since $2\leq i \leq n-2$, then we can write $y_i=[y_{0},y_{i-1}]$. By using the Leibniz identity and the hypothesis of the corollary we have
$[y_{1},y_i]=-[y_{i},y_1].$
\item{ Let us see $[y_{j},y_i]=-[y_i,y_j] $ with $ 2\leq i,j \leq n-2.$}\\
It is enough to apply the Corollary \ref{Lema1}.
\end{itemize}
\end{dem}

We are now in a position to prove that it is impossible to obtain a Leibniz non-Lie algebra of maximum length
 from the generalization of naturally graded quasi-filiform Lie algebras by using natural graduation.

\begin{thm}\label{Lieext}
Every quasi-filiform Leibniz algebra obtained by natural graduation extension of naturally graded quasi-filiform Lie algebra is a Lie algebra.
\end{thm}
\begin{dem}
Since all of the above algebras are quasi-filiform, $[x_0,x_i]=x_{i+1},$ where $ 1\leq i \leq n-3,$ is satisfied. Then, when we extend them using natural
 graduations, it can be considered a new basis with the following generators: $$x_s=x_0+\sum_{i=1}^{n-1}a_ix_i;\qquad x_t=x_1+\sum_{j=0, j\neq 1}^{n-1}b_jx_j ; \qquad x_u=x_{n-1}+\sum_{k=0}^{n-2}c_kx_k,$$
 and the following products are obtained:
\begin{align*}
[x_s,x_s]&=(*)x_3+ \dots+(*)x_{n-1},\\
[x_t,x_t]&=(*)x_3+ \dots+(*)x_{n-1},\\
[x_u,x_u]&=(*)x_3+ \dots+(*)x_{n-1},
\end{align*}

Moreover, we will only work with the products $[x_s,x_t]$, $[x_s,x_u]$ and $[x_t,x_u]$, because $[x_i,x_j]=-[x_j,x_i] \quad \forall x_i,x_j \in L_1$.

The generalization of that naturally graded quasi-filiform Lie algebras always verifies:
\begin{align*}
[x_0,x_i]&=x_{i+1}+(*)x_{i+2}+ \dots+(*)x_{n-1},\\
[x_i,x_0]&=-x_{i+1}+(\diamond)x_{i+2}+ \dots +(\diamond)x_{n-1}.
\end{align*}

Then, if we take the products: $$
[\underbrace{[x_s,\dots,[x_s,[x_s,}_{\hbox{i-times}}x_t]]]=(1-a_1b_0)x_{i+1}+(*)x_{i+2}+
\dots+x_{n-1} $$ or
$$[\underbrace{[x_s,\dots,[x_s,[x_s,}_{\hbox{i-times}}x_u]]]=(c_1-c_0a_1)x_{i+1}+(*)x_{i+2}+
\dots+(*)x_{n-1},$$
or
$$[\underbrace{[x_s,\dots,[x_s}_{\hbox{(i-1)-times}},[x_u,x_t]]]=(c_0-c_1b_0)x_{i+1}+(*)x_{i+2}+ \dots+(*)x_{n-1}$$
\noindent with $2\leq i \leq n-4$, the basis $b_1$ can be constructed:
$$b_1=\{y_0=x_s,\ y_1=x_t,\ y_i=[y_0,y_{i-1}],\ 2 \leq i \leq n-3,\ y_{n-2}=?, \ y_{n-1}=?\}$$

\noindent or $b_2$ if there are three generators:
$$b_2=\left\{y_0=x_s,\ y_1=x_u,\ y_{n-1}=x_t, \ y_i=[y_0,y_{i-1}],\hbox{ with } \ 2 \leq i \leq n-3, \ y_{n-2}=?\right\}.$$

From its very definition and taking into account property $[x_i,x_j]=-[x_j,x_i] $ with $x_i,\ x_j \in L_1$ it can be affirmed
$[y_i,y_0]=-[y_0,y_i]=-y_{i+1}$ with $ 2 \leq i \leq n-3$.
Hence, it is enough to choose $y_{n-2}$ and $y_{n-1}$ (when $y_{n-1}$ is not the generator) and to prove the corresponding corollary in each case.

Let us study two main algebras which allow us to prove the other ones.

\

{\it The algebra $\widetilde{L_{n-1}} \oplus \mathbb{C}$. (Analogously for $\widetilde{\mathrm{Q}_{n-1}} \oplus \mathbb{C}$ ).}

\

If we consider the natural graduation $L_{1}=<x_0,x_1,x_{n-1}>$ and $L_{i}=<x_i> \hbox{ for } 2 \leq i \leq n-2$, then
the generalization of $L_{n-1}$ is defined by:

$$\begin{cases}
[x_0,x_i]=x_{i+1}+(*)x_{i+2}+\dots+(*)x_{n-2},\\
[x_i,x_0]=-x_{i+1}+(\diamond)x_{i+2}+ \dots+(\diamond)x_{n-2},\\
[x_i,x_j]=(*)x_{i+j+1}+\dots+(*)x_{n-2} \quad (i,j)\neq (0,j),(i,0).
\end{cases}$$

Three cases can be considered by using the maximum length study:\\

\textbf{Case 1: If $1-a_1b_0 \neq 0.$}

 Let us take the basis $b_1 $ with $y_{n-2}=[y_0,y_{n-3}$, $y_{n-1}=x_u$ and its graduation of maximum length:
 $V_{k_s}\oplus V_{k_t}\oplus V_{k_t+k_s}\oplus V_{k_t+2k_s}\oplus \dots \oplus V_{k_t+(n-3)k_s}\oplus V_{k_u}$.

 By the corollary \ref{Lema1} and the similar arguments to Theorem \ref{th:tipo I}, we can conclude that the algebra is a Lie algebra.

\textbf{Case 2: If $c_0a_1-c_1 \neq 0,$ $1-a_1b_0 =0.$}

 Let us take the basis $b_2$ with $y_{n-2}=[y_0,y_{n-3}]$. In order to prove that the algebra is a Lie algebra, it is sufficient to show that the hypotheses of Corollary \ref{Lema1} are true.
That is obtained by a similar argument to the previous one.

  \textbf{Case 3: If $c_0-c_1b_0 \neq 0,$ $c_0a_1-c_1=1-a_1b_0 =0.$}

  This case is impossible.

\

 {\it The algebra $\widetilde{\tau_{(n,n-3)}}.$ (Analogously for $\widetilde{\tau_{(n,n-4)}}).$}

 \

       The natural graduation of $\tau_{(n,n-4)}$ is formed by the
     subset $L_1=<x_0,x_1>$, $L_i=<x_i> $ with $2 \leq i \leq
     n-4$, $L_{n-3}=<x_{n-3},x_{n-1}>$ and $L_{n-2}=<x_{n-2}>$.
     Its natural generalization is defined by the following
     products::
     $$\left\{\begin{array}{ll}
     [x_0,x_i]=x_{i+1}+(*)x_{i+2}+ \dots+(*)x_{n-1},& 1 \leq i \leq n-5,\\{}
     [x_0,x_{n-4}]=x_{n-3}+(*)x_{n-2},&\\{}
     [x_0,x_{n-3}]=x_{n-2},&\\{}
     [x_i,x_0]=-x_{i+1}+(*)x_{i+2}+ \dots+(*)x_{n-1}, &1 \leq i \leq n-5,\\{}
     [x_{n-4},x_{0}]=-x_{n-3}+(*)x_{n-2},&\\{}
     [x_{n-3},x_{0}]=-x_{n-2},&\\{}
     [x_{n-1},x_{1}]=\frac{n-4}{2}x_{n-2},&\\{}
     [x_1,x_{n-1}]=- \frac{n-4}{2}x_{n-2},&\\{}
     [x_i,x_{n-3-i}]=(-1)^{i-1}(x_{n-3}+x_{n-1})+(*)x_{n-2}, & 1 \leq i \leq \frac{n-4}{2},\\{}
     [x_{n-3-i},x_{i}]=(-1)^{i}(x_{n-3}+x_{n-1})+(*)x_{n-2}, & 1 \leq i \leq \frac{n-4}{2},\\{}
     [x_i,x_{n-2-i}]=(-1)^{i-1}\frac{n-2-2i}{2}x_{n-2},& 1 \leq i \leq \frac{n-4}{2},\\{}
     [x_{n-2-i},x_{i}]=(-1)^{i}\frac{n-2-2i}{2}x_{n-2}, & 1 \leq i \leq \frac{n-4}{2},\\{}
     [x_i,x_j]=(*)x_{i+j+1}+\dots+(*)x_{n-2}, & \hbox{in other cases}.
     \end{array}\right.$$

     We consider the basis $b_1$. It follows from the independence of the generator of the basis that $1-a_1b_0\neq 0$.

     We will need to choose $y_{n-3}$, $y_{n-2}$ and $y_{n-1}$ and to prove the hypotheses of Corollary \ref{Lema1}.
     We can consider two cases to choose these vectors.

    \textbf{Case 1: If $1+\frac{n-2}{2} a_1 \neq 0, $} then $y_{n-2}=[y_0,y_{n-3}]$ and $y_{n-1}=[y_1,y_{n-4}]$ can be constructed.

     \textbf{Case 2: If $1+a_1\frac{n-2}{2}=0, $} then $y_{n-1}=[y_1,y_{n-4}]$ and $y_{n-2}=[y_0,y_{n-1}]$ can be constructed.

     We can use the same argument applied in the previous cases to prove that the hypotheses of Corollary \ref{Lema1} are true,
     hence the algebra is a Lie algebra.

 \end{dem}

\section {$4$-dimensional and $5$-dimensional quasi-filiform Leibniz algebras of maximum length}
In order to classify quasi-filiform Leibniz algebras of maximum length the dimension $n$ must be $\geq 6$. Consequently, in this section we will show the result of $n=4 $ and $n=5$. We follow the same method as in the general case.

\begin{thm}\label{n=4}  Let $\ll$ be a $4$-dimensional Leibniz algebra of maximum length, then $\ll$ is isomorphic to one algebra of the following pairwise non isomorphic families:
$$\begin{array}{ll}
N^{1,\alpha}: & N^{2}:\\[2mm]
\begin{cases}
[y_1,y_1]=y_2,\\
[y_3,y_1]=y_4,\\
[y_1,y_3]=\alpha y_4, \quad \alpha \in \mathbb{C}.
\end{cases} & \begin{cases}
               [y_1,y_1]=y_2,\\
               [y_1,y_3]=y_4.
             \end{cases}
   \end{array}$$
   \end{thm}

\begin{thm}\label{n=5}
   Let $\ll$ be a 5-dimensional Leibniz algebra of maximum length, then $\ll$ is isomorphic to one algebra of the following pairwise non isomorphic families:\\
   $$\begin{array}{lll}

  M^{1,0}:& M^{2,\lambda}:&
     M^{3,0}:\\
     
  \begin{cases}
            [y_1,y_1]=y_5,\\
            [y_{4},y_1]=y_2,\\
            [y_2,y_1]=y_{3}.\\
           \end{cases} & \begin{cases}
                 [y_i,y_1]=y_{i+1}, &1\leq i \leq 2\\
                 [y_{4},y_1]= y_5,\\
                 [y_1,y_{4}]=\lambda y_{5}, &\lambda \in \mathbb{C}.
                \end{cases}&

\begin{cases}
                   [y_1,y_1]=y_2\\
                   [y_i,y_1]=y_{i+1} &3\leq i \leq 4\\
                   [y_1,y_4]=-y_{5}
                 \end{cases}

   \end{array}$$
\end{thm}

\

Let us remark that the previous theorems complete the study of
Leibniz algebras of maximum length with nilindex up to $n-2$ with $n\geq 4$.


\textbf{Acknowledgments.} \emph{This work is
supported in part by the PAI, FQM143 of Junta de
Andaluc\'{\i}a (Spain). B.A. Omirov was supported by a grant of
NATO-Reintegration ref. CBP.EAP.RIG. 983169 and he would like to thank of
the Universidad de Sevilla for their hospitality.}

\

{\sc Luisa M. Camacho, Elisa M. Ca\~{n}ete, Jos\'{e} R. G\'{o}mez.}  Dpto. Matem\'{a}tica Aplicada I.
Universidad de Sevilla. Avda. Reina Mercedes, s/n. 41012 Sevilla.
(Spain), e-mail: \emph{lcamacho@us.es}, \emph{elisacamol@us.es}, \emph{jrgomez@us.es}

\

{\sc Bakhrom A. Omirov.} Institute of Mathematics and Information Technologues, Uzbekistan
Academy of Science, F. Hodjaev str. 29, 100125, Tashkent
(Uzbekistan), e-mail: \emph{omirovb@mail.ru}

\end{document}